\documentclass[11pt,a4paper]{article}
\usepackage{amsmath}
\usepackage{amssymb}
\usepackage{amsthm}
\usepackage{amsfonts}
\usepackage{enumerate}
\usepackage{verbatim}
\usepackage{hyperref}
\usepackage{breakurl,url}
\usepackage{bm}
\usepackage{color}
\usepackage{mathtools} %for := \coloneqq, for parentheses
\usepackage{pstricks}
\usepackage{pst-node}
\usepackage{pst-plot}
\usepackage{algorithm}
\usepackage{algorithmic}
\usepackage{pst-poly}
\usepackage{caption,subcaption} %subfigures
\usepackage[margin=1.3in]{geometry}
\newcommand{\mmid}[0]{;\,}		%

\newcommand{\omace}[1]{\mbox{$\overline{#1}$}}	%upper limit of int matrix
\newcommand{\umace}[1]{\mbox{$\underline{#1}$}}  %lower limit of int matrix
\DeclareMathOperator{\diag}{diag}	%diagonal matrix
\DeclareMathOperator{\sgn}{sgn}	%sign
\DeclareMathOperator{\vek}{vec} %vectorization of a matrix

\newcommand{\R}[0]{{\mathbb{R}}}
\newcommand{\C}[0]{{\mathbb{C}}}
\def\clqq{``}
\def\crqq{''}
\def\quo#1{\clqq{}#1\crqq{}}  % en. quotes

\newtheorem{theorem}{Theorem}

\newtheorem{corollary}{Corollary}

\theoremstyle{definition}

\newtheorem{example}{Example}
\newtheorem{remark}{Remark}
\begin{document}

\title{
%\textcolor{red} {Necessary and sufficient conditions for unique solvability of absolute value equations: An overview, extensions and perspectives}
%\textcolor{blue} {Investigating Unique Solvability of Absolute Value Equations: An Overview, Extensions, and Future Directions}
Characterization of Unique Solvability of Absolute Value Equations: An Overview, Extensions, and Future Directions}

\author{
Shubham Kumar\footnote{
Department of Mathematics, PDPM-Indian Institute of Information Technology, Design and Manufacturing Jabalpur, M.P. India, e-mail: \texttt{shub.srma@gmail.com}}
 \and Deepmala\footnote{
Department of Mathematics, Faculty of Natural Sciences, PDPM-Indian Institute of Information Technology, Design and Manufacturing Jabalpur, M.P. India, e-mail: \texttt{dmrai23@gmail.com}}
 \and  Milan Hlad\'{i}k\footnote{
Charles University, Faculty  of  Mathematics  and  Physics,
Department of Applied Mathematics,
Malostransk\'e n\'am.~25, 11800, Prague, Czech Republic,
e-mail: \texttt{hladik@kam.mff.cuni.cz}}
\and Hossein Moosaei\footnote{
Department of Informatics, Faculty of Science, Jan Evangelista Purkyně University, \'{U}st\'{i} nad Labem, Czech Republic, 
e-mail: \texttt{Hossein.Moosaei@ujep.cz}}
}

\date{\today}
\maketitle

\begin{abstract}
This paper provides an overview of the necessary and sufficient conditions for guaranteeing the unique solvability of absolute value equations. In addition to discussing the basic form of these equations, we also address several generalizations, including generalized absolute value equations and matrix absolute value equations. Our survey encompasses known results as well as novel characterizations proposed in this study.
\end{abstract}

\bigskip
%\paragraph{MSC Codes.}   90C30, 90C26, 90C33, 65G40,  15A06

%%%%%%%%%%%%%%%%%%%%%%%%%%%%%%%%%%%%%%%%%%%%%%%%%%%%%%%%%%%%%%%
% INTRODUCTION
%%%%%%%%%%%%%%%%%%%%%%%%%%%%%%%%%%%%%%%%%%%%%%%%%%%%%%%%%%%%%%%
\section{Introduction}

The study of absolute value systems of equations in various forms has garnered significant attention in both theoretical and algorithmic research over the past two decades \cite{HlaMos2022u,ketabchi2012minimum,mangasarian2009generalized,moosaei2021optimal,raayatpanah2019absolute,zamani2021new}.  This paper aims to provide an overview on one particular but essential topic -- the unique solvability characterization of such systems. To facilitate our discussion, we begin by introducing the required notation and key concepts.

%%%%%
\paragraph{Notation.}
We use $A_{*i}$ to denote $i$-th column of a matrix~$A$ and $\diag(v)$ for the diagonal matrix with entries given by vector $v\in\R^n$. The identity matrix of size $n\times n$ is denoted by~$I_n$. The sign of a real $r$ is $\sgn(r)=1$ if $r\geq0$ and $\sgn(r)=-1$ otherwise. 
For vectors and matrices, we understand the inequalities, sign function and absolute values entrywise.

%%%%%
\paragraph{Definitions.}
A matrix $R$ is \emph{a column representative} of a pair of matrices $\{C,D\}$ if $R_{*i}\in\{C_{*i},D_{*i}\}$ for every~$i$. Next, a pair of matrices $\{C,D\}$ has \emph{the column $\mathcal{W}$-property} if all its column representative matrices have the positive determinant, or all of them have the negative determinant. Similarly, we define \emph{the row $\mathcal{W}$-property.}

Let $\umace{A},\omace{A}\in\R^{m\times n}$ such that $\umace{A}\leq\omace{A}$. \emph{The interval matrix} is defined as the set of matrices
\begin{align*}
[\umace{A},\omace{A}]=\{A\mmid \umace{A}\leq A\leq \omace{A} \}.
\end{align*}
We say that $[\umace{A},\omace{A}]$ is regular if every $A\in[\umace{A},\omace{A}]$ is nonsingular. A survey on regularity of interval matrices is provided in Rohn~\cite{Roh2009}.

A matrix $A\in\R^{n\times n}$ is a \emph{P-matrix} if every principal minor (i.e., the determinant of a principal submatrix) is positive.

%%%%%
\paragraph{Preliminaries.}
The standard system of \emph{absolute value equations (AVE)} has the form of
\begin{align}\label{ave}
Ax+|x| = b,
\end{align}
where $A\in \mathbb{R}^{n\times n}$, $b \in \mathbb{R}^n$ are given and $x\in\R^n$ is the unknown. 
Due to the absolute value, checking its solvability is NP-hard~\cite{Man2007}. Also, other related problems are intractable~\cite{Pro2009}.

It is well-known \cite{Man2007,ManMey2006} that AVE is equivalent to the linear complementarity problem (LCP), which is formulated as follows:
\begin{align}\label{lcp}
y=Mz+q,\ y^Tz=0,\ y,z\geq0,
\end{align}
where $M\in\R^{n\times n}$ and $q\in\R^n$ is an input and $y,z\in\R^n$ are variables. The LCP has a unique solution for every $q\in\R^n$ if and only if the constraint matrix $M\in \mathbb{R}^{n\times n}$ is a  P-matrix; see, e.g., Cottle et al.~\cite{CotPan2009}. By Coxson~\cite{Cox1994}, checking P-matrix property is a co-NP-hard problem, which makes this issue hard and also challenging. 

Besides AVE, \emph{the generalized absolute value equations (GAVE)} 
have also been extensively studied. The GAVE can be represented as follows:
\begin{align}\label{gave}
Ax+B|x| = b,
\end{align}
where $A,~ B\in \mathbb{R}^{n\times n}$ are matrices, and $x,~ b \in \mathbb{R}^n$ are vectors. Similar to the equivalence of AVE with the LCP, it has been shown by Mezzadri~\cite{Mez2020} that the GAVE is equivalent to a specific variant of the LCP known as the horizontal LCP.

In general, the solution set of GAVE is neither convex nor connected. However, in any orthant, it forms a convex polyhedron. This is easy to see since the orthant given by the sign vector $s\in\{\pm1\}^n$ is characterized by $\diag(s)x\geq0$ and the absolute value takes the form $|x|=\diag(s)x$. Therefore, the solution set restricted to the orthant reads as
\begin{align}\label{gaveOrth}
(A+B\diag(s))x=b,\ \ \diag(s)x\geq0.
\end{align}

Our aim is to characterize the  unique solvability of AVE and GAVE for each right-hand side vector $b\in\R^n$. In view of the equivalence between AVE and LCP, these problems are intractable. So the conditions we present in the next sections are hard to check and often of exponential manner. That is why researchers have also been interested in computationally cheap sufficient conditions (see, e.g., \cite{Hla2023a,HlaMoo2023a,ManMey2006,RohHoo2014,WuShe2021,WuGuo2016}), but our focus is merely on the complete characterization. We also do not discuss sufficient conditions for unsolvability \cite{Hla2018b,ManMey2006} and sufficient conditions for existence of $2^n$ solutions \cite{Hla2018b,ManMey2006}, which is the largest finite number of solutions that may exist for AVE or GAVE.

%%%%%%%%%%%%%%%%%%%%%%%%%%%%%%%%%%%%%%%%%%%%%%%%%%%%%%%%%%%%%%%
\section{Absolute value equations (AVE)}\label{sAve}

A convenient way to characterize the unique solvability of AVE is through the use of interval matrices. Theorem~\ref{thmAveUnSol}, which is explicitly stated by Wu and Li~\cite{WuLi2018}, provides such a characterization. However, partial results regarding this characterization have already been presented in previous works, such as \cite{Roh2004b, Roh2012c, Rum2003, ZhaWei2009}.

\begin{theorem}\label{thmAveUnSol}
The AVE has a unique solution for each $b\in\R^n$ if and only if $[A-I_n,A+I_n]$ is regular.
\end{theorem}

Regularity of interval matrices was surveyed by Rohn~\cite{Roh2009}, who presented 40 necessary and sufficient conditions for regularity of interval matrices. Essentially, all of these conditions can be used to characterize the regularity of $[A-I_n,A+I_n]$ and, consequently, the unique solvability of AVE. However, for the sake of illustration, we will only list a selection of them.

\begin{corollary}
The following conditions are equivalent:
\begin{enumerate}[(i)]
%\item 
%The AVE has a unique solution for each $b\in\R^n$;
\item
$[A-I_n,A+I_n]$ is regular;
\item
the system $|Ax|\leq|x|$ has only the trivial solution $x=0$;
\item
for each $s\in\{\pm1\}^n$, the linear system
$$
(A-\diag(s))x^1-(A+\diag(s))x^2=s,\ x^1,x^2\geq0
$$
is solvable;
\item
$\det(A+\diag(s))$ is constantly positive or negative for each $s\in\{\pm1\}^n$;
\item
each matrix in the form 
$A+t\diag(s)$
is nonsingular, where $s\in\{\pm1\}^n$ and $t\in[0,1]$;
\item
%each matrix in the form $A+\diag(z)$ is nonsingular, where $z_i\in[-1,1]$ for some $i\in\{1,\ldots,n\}$ and $|z_j|=1$ for $j\not=i$.
each matrix in the form 
$$A+\diag(z^{(i)}),\quad i=1,\ldots,n,$$
is nonsingular, where $z^{(i)}_i\in[-1,1]$ and $|z^{(i)}_j|=1$ for $j\not=i$;
\item
%$A$ is nonsingular and for each $s\in\{\pm1\}^n$ and each real eigenvalue $\lambda$ of $A^{-1}\diag(s)$ we have $|\lambda|<1$.
for each $s\in\{\pm1\}^n$ and each real eigenvalue $\lambda$ of $\diag(s)A$ we have $|\lambda|>1$.
\end{enumerate}
\end{corollary}

\begin{proof}
The conditions follow from Rohn \cite{Roh2009}, Theorem~4.1., conditions (i), (ii), (xxv), (xxxii), (xxxviii), (xli), and (xxxiv), respectively.
\end{proof}

Based on the reduction of AVE to LCP, we immediately have the following characterization \cite{ManMey2006,Mez2020,WuLi2018}.

\begin{theorem}\label{thmAveUnSolPmat}
The AVE has a unique solution for each $b\in\R^n$ if and only if $(A+I_n)^{-1}(A-I_n)$ is a P-matrix.
\end{theorem}

As a direct consequence of Theorem~\ref{thmAveUnSol}, we can deduce under which transformation of the constraint matrix the unique solvability is preserved.

\begin{corollary}\label{corAveUnTrans}
The following conditions are equivalent:
\begin{enumerate}[(i)]
\item 
The system $Ax+|x|=b$ has a unique solution for each $b\in\R^n$;
\item
the system $A^Tx+|x|=b$ has a unique solution for each $b\in\R^n$;
\item
the system $\alpha Ax+|x|=b$ has a unique solution for each $b\in\R^n$ and each $|\alpha|\geq1$.
\end{enumerate}
\end{corollary}

The following example demonstrates that while certain transformations (such as those mentioned in Corollary~\ref{corAveUnTrans}) can preserve the unique solvability of the AVE, other transformations (such as squaring the matrix) may result in a loss of unique solvability. This emphasizes the necessity of analyzing the properties of transformed matrices carefully to ensure the preservation of unique solutions in practical applications of the AVE.

\begin{example}
The unique solvability is not preserved under the square operation on matrix~$A$. For example, let
$$
A=\begin{pmatrix}0&2\\-1&1\end{pmatrix}.
$$
Then $[A-I_2,A+I_2]$ is not regular since it contains the singular matrix 
$
\begin{psmallmatrix}-1&2\\-1&2\end{psmallmatrix}.
$
On the other hand, the interval matrix 
$$
[A^2-I_2,A^2+I_2]=\begin{pmatrix}[-3,-1]&2\\-1&[-2,0]\end{pmatrix}
$$ 
is regular.

Also the converse situation may occur. Let
$$
A=\begin{pmatrix} -1 & 1 & -2\\ -3 & 3 & -3\\ -1 & 3 & -3\end{pmatrix}.
$$
Now, $[A-I_3,A+I_3]$ is regular. However, the interval matrix $[A^2-I_3,A^2+I_3]$ 
is not regular since it contains the singular matrix 
$$
\begin{pmatrix} 1 & -4 & 5\\ -3 & -3.75 & 6\\ -5 & -1 & 3\end{pmatrix}\in
\begin{pmatrix} [-1,1] & -4 & 5\\ -3 & -[2,4] & 6\\ -5 & -1 & [1,3]\end{pmatrix}
=[A^2-I_3,A^2+I_3].
$$ 
\end{example}

\begin{remark}
Characterization of unique solvability for particular matrix classes (such as symmetric, nonnegative or positive definite matrices) was addressed, e.g., in \cite{Hla2023a,HlaMoo2023a}. However, in this survey paper, we focus merely on general matrices.
\end{remark}

%%%%%%%%%%%%%%%%%%%%%%%%%%%%%%%%%%%%%%%%%%%%%%%%%%%%%%%%%%%%%%%
\section{Generalized absolute value equations (GAVE)}\label{sGave}

Theorem~\ref{thmAveUnSol} was extended for the case of GAVE by Wu and Shen~\cite{WuShe2021}.

\begin{theorem}\label{thmGaveUnSol}
The GAVE has a unique solution for any $b\in\R^n$ if and only if matrix $A + BD$ is nonsingular for every $D\in[-I_n,I_n]$.
\end{theorem}

Notice that the set of matrices
\begin{align}\label{mnaMcABD}
\{A + BD\mmid D\in[-I_n,I_n]\}
\end{align}
does not form a standard interval matrix. This type of sets is usually called a linear interval parametric matrix, and they are much more difficult to deal with; see, e.g., \cite{May2017,Ska2018}. Nevertheless, due to its special structure, we can derive an explicit characterization.

\begin{theorem}\label{thmGaveCharT}
Matrix $A + BD$ is nonsingular for every $D\in[-I_n,I_n]$ if and only if the system
\begin{align}\label{ineqCharRegM}
|A^Tx| \leq  |B^Tx| 
\end{align}
has only the trivial solution $x=0$.
\end{theorem}

\begin{proof}
Equivalently, we consider for nonsingularity the set of matrices $A^T + DB^T$, $D\in[-I_n,I_n]$. Since each interval parameter $d_{ii}$ influences only one row of the resulting parametric matrix, we have by Popova~\cite{Pop2009} that the set of matrices is nonsingular if and only if the system
\begin{align*}
|A^Tx| \leq \sum_{i=1}^n e_i|(B^T)_{i*}x|
\end{align*}
has only the trivial solution $x=0$. The right-hand side of the system equivalently reads as 
$|B^Tx|$. 
\end{proof}

Another type of characterization follows from a transformation of the set \eqref{mnaMcABD} to a standard (albeit double-sized) interval matrix.

\begin{theorem}\label{thmGaveRegABD}
Matrix $A + BD$ is nonsingular for every $D\in[-I_n,I_n]$ if and only if the interval matrix
\begin{align}\label{intMcDouble}
\begin{pmatrix}A&B\\{}[-I_n,I_n]&I_n\end{pmatrix}
\end{align}
is regular.
\end{theorem}

\begin{proof}
Let $D\in[-I_n,I_n]$ be arbitrary. By elementary operations
\begin{align*}
\begin{pmatrix}A&B\\D&I_n\end{pmatrix}
\sim\begin{pmatrix}B&A\\I_n&D\end{pmatrix}
\sim\begin{pmatrix}0&A-BD\\I_n&D\end{pmatrix}
\end{align*}
we see the equivalence.
\end{proof}

Regularity of interval matrix \eqref{intMcDouble} can again be characterized by means of 40 conditions by Rohn~\cite{Roh2009}; we adapt a selection of them.

\begin{corollary}\label{corIntMcDouble}
The following conditions are equivalent:
\begin{enumerate}[(i)]
%\item 
%The AVE has a unique solution for each $b\in\R^n$;
\item\label{itCorIntMcDouble1}
Interval matrix \eqref{intMcDouble} is regular;
\item\label{itCorIntMcDouble2}
the system $Ax+By=0$, $|y|\leq|x|$ has only the trivial solution $x=y=0$;
\item\label{itCorIntMcDouble3}
for each $s,t\in\{\pm1\}^n$, the linear system
$$
A(x^1-x^2)+B\diag(s)(x^1+x^2)=t-B\diag(s),\ x^1,x^2\geq0
$$
is solvable;
\item\label{itCorIntMcDouble4}
$\det(A+B\diag(s))$ is constantly positive or negative for each $s\in\{\pm1\}^n$;
\item\label{itCorIntMcDouble5}
each matrix in the form 
$A+tB\diag(s)$
is nonsingular, where $s\in\{\pm1\}^n$ and $t\in[0,1]$;
\item\label{itCorIntMcDouble6}
%each matrix in the form $A+\diag(z)$ is nonsingular, where $z_i\in[-1,1]$ for some $i\in\{1,\ldots,n\}$ and $|z_j|=1$ for $j\not=i$.
each matrix in the form 
$$A+B\diag(z^{(i)}),\quad i=1,\ldots,n,$$
is nonsingular, where $z^{(i)}_i\in[-1,1]$ and $|z^{(i)}_j|=1$ for $j\not=i$;
\item\label{itCorIntMcDouble6b}
%$A$ is nonsingular and for each $s\in\{\pm1\}^n$ and each real eigenvalue $\lambda$ of $A^{-1}\diag(s)$ we have $|\lambda|<1$.
$A$ is nonsingular and  for each $s\in\{\pm1\}^n$ and each real eigenvalue $\lambda$ of $\diag(s)A^{-1}B$ we have $|\lambda|<1$;
\item\label{itCorIntMcDouble7}
the map $x\mapsto Ax+B|x|$ is a bijection.
\end{enumerate}
\end{corollary}

\begin{proof}
\quo{Condition \eqref{itCorIntMcDouble2}.}
By Rohn \cite[cond.~(ii)]{Roh2009}, \eqref{intMcDouble} is regular if and only if the system $|Ax+By|\leq0$, $|y|\leq|x|$ has only the trivial solution.

\quo{Condition \eqref{itCorIntMcDouble3}.}
By Rohn \cite[cond.~(xxv)]{Roh2009}, \eqref{intMcDouble} is regular if and only if for each $s,t\in\{\pm1\}^n$, the linear system
$$
A(x^1-x^2)+B(y^1-y^2)=t,\ 
-\diag(s)(x^1+x^2)+(y^1-y^2)=s,\ 
x^1,x^2,y^1,y^2\geq0
$$
is solvable. Substituting for $(y^1-y^2)$ from the second equation, we obtain the required system.

\quo{Condition \eqref{itCorIntMcDouble4}.}
This follows from Rohn \cite[cond.~(xxxii)]{Roh2009}, or directly by linearity of the determinant in each column.

\quo{Condition \eqref{itCorIntMcDouble5}.}
This follows from Rohn \cite[cond.~(xxxviii)]{Roh2009}.

\quo{Condition \eqref{itCorIntMcDouble6}.}
This follows from Rohn \cite[cond.~(xli)]{Roh2009}.

\quo{Condition \eqref{itCorIntMcDouble6b}.}
This is a consequence of subsequent Theorem~\ref{thmGaveSpecUnSolv}\eqref{itThmGaveSpecUnSolv5} since GAVE is under the assumption equivalent to $x+A^{-1}B|x|=A^{-1}b$.

\quo{Condition \eqref{itCorIntMcDouble7}.}
This is a reformulation of unique solvability of GAVE $Ax+B|x|=b$ for every $b\in\R^n$.
\end{proof}

Based on the equivalence of GAVE and the horizontal LCP, Mezzadri~\cite{Mez2020} deduced another equivalent conditions on unique solvability. 

\begin{theorem}
The following conditions are equivalent:
\begin{enumerate}[(i)]
\item 
The GAVE has a unique solution for each $b\in\R^n$;
\item
$\{A-B,A+B\}$ has the column $\mathcal{W}$-property;
\item
$A+B$ is nonsingular and $\{I_n,(A+B)^{-1}(A-B)\}$ has the column $\mathcal{W}$-property;
\item
$A+B$ is nonsingular and $(A+B)^{-1}(A-B)$ is a P-matrix.
\end{enumerate}
\end{theorem}

What are the transformations of $A$ and $B$ preserving the unique solvability?

\begin{corollary}\label{corGaveUnSolvTrans}
The following conditions are equivalent:
\begin{enumerate}[(i)]
\item 
The system $Ax+B|x|=b$ has a unique solution for each $b\in\R^n$;
%\item
%for some $D:|D|=I_n$, the system $DAx+B|x|=b$ has a unique solution for each $b\in\R^n$;
%\item
%for every $D:|D|=I_n$, the system $DAx+B|x|=b$ has a unique solution for each $b\in\R^n$;
\item\label{itCorGaveUnSolvTrans1}
for some $D:|D|=I_n$, the system $ADx+B|x|=b$ has a unique solution for each $b\in\R^n$;
\item\label{itCorGaveUnSolvTrans2}
for every $D:|D|=I_n$, the system $ADx+B|x|=b$ has a unique solution for each $b\in\R^n$;
\item\label{itCorGaveUnSolvTrans3}
for some $D:|D|=I_n$, the system $Ax+BD|x|=b$ has a unique solution for each $b\in\R^n$;
\item\label{itCorGaveUnSolvTrans4}
for every $D:|D|=I_n$, the system $Ax+BD|x|=b$ has a unique solution for each $b\in\R^n$;
\item\label{itCorGaveUnSolvTrans5}
for every $D\in[-I_n,I_n]$, the system $Ax+BD|x|=b$ has a unique solution for each $b\in\R^n$;
%\item
%the system $A^Tx+B^T|x|=b$ has a unique solution for each $b\in\R^n$;
\item\label{itCorGaveUnSolvTrans6}
the system $\alpha Ax+B|x|=b$ has a unique solution for each $b\in\R^n$ and each $|\alpha|\geq1$.
\end{enumerate}
\end{corollary}

\begin{proof}
Conditions \eqref{itCorGaveUnSolvTrans1}--\eqref{itCorGaveUnSolvTrans2} follow from the fact that the class $A+BD'$ is nonsingular for each $D'\in[-I_n,I_n]$ if and only if $AD+BD'$ is nonsingular for each $D'\in[-I_n,I_n]$.

Similarly conditions \eqref{itCorGaveUnSolvTrans3}--\eqref{itCorGaveUnSolvTrans5} follow from the fact that the class $A+BD'$ is nonsingular for each $D'\in[-I_n,I_n]$ if and only if $A+BDD'$ is nonsingular for each $D'\in[-I_n,I_n]$.

Condition \eqref{itCorGaveUnSolvTrans6} is obvious in view of Theorem~\ref{thmGaveUnSol}.
\end{proof}

\begin{example}
The statements of Corollary~\ref{corGaveUnSolvTrans} are not valid for the transformation 
$Ax+DB|x|=b$, where $D$ is such that $|D|=I_n$. To see it, consider
$$
A=\begin{pmatrix}-1&1\\1&1\end{pmatrix},\ \ 
B=\begin{pmatrix}1&0\\1&0\end{pmatrix},\ \ 
D=\begin{pmatrix}-1&0\\0&1\end{pmatrix}.
$$
Then matrix $A+BD'$ is nonsingular for each $D'\in[-I_n,I_n]$ (its determinant is negative). Nevertheless, this is not true for matrices of type $A+DBD'$, $D'\in[-I_n,I_n]$, since they cover the singular matrix $\begin{psmallmatrix}0&1\\0&1\end{psmallmatrix}$.

For this reason, the unique solvability is not preserved under matrix transposition. Thus, $Ax+B|x|=b$ can be uniquely solvable for each $b\in\R^n$, but $A^Tx+B^T|x|=b$ need not (in contrast to AVE or the special case $A=I_n$ discussed below in Corollary~\ref{corGaveAspec}).
\end{example}

Eventually, we present several equivalent conditions related to the type of solvability and the number of solutions.

\begin{theorem}\label{thmGaveSolvVar}
The following conditions are equivalent:
\begin{enumerate}[(i)]
%\item 
%The AVE has a unique solution for each $b\in\R^n$;
\item\label{itThmGaveSolvVar1}
The system $Ax+B|x|=b$ has a unique solution for each $b\in\R^n$;
\item\label{itThmGaveSolvVar2}
The system $Ax+B|x|=b$ has at most one solution for each $b\in\R^n$;
\item\label{itThmGaveSolvVar3}
The system $Ax+BD|x|=b$ has at least one solution for each $b\in\R^n$ and each $D:|D|=I_n$.
\end{enumerate}
\end{theorem}

\begin{proof}
\quo{Condition \eqref{itThmGaveSolvVar2}.}
This is a consequence of subsequent Theorem~\ref{thmGaveSpecUnSolv}\eqref{itThmGaveSpecUnSolv1}. We simply rewrite $Ax+B|x|=b$ as $x+A^{-1}B|x|=A^{-1}b$. Notice that $A$ must be nonsingular. Otherwise, there is $x^*\not=0$ such that $Ax^*=0$. In this case we put $b\coloneqq B|x^*|$, and the system $Ax+B|x|=b$ has at least two solutions, $x=\pm x^*$.

\quo{Condition \eqref{itThmGaveSolvVar3}.}
If condition \eqref{itThmGaveSolvVar1} is not true, then by Theorem~\ref{thmGaveCharT} there exists $y^*\not=0$ such that $|A^Ty^*| \leq  |B^Ty^*|$. Define $b\coloneqq -y^*$ and $D\coloneqq\diag(\sgn(B^Ty^*))$. Then the system
\begin{align*}
DB^Ty \geq -\diag(s)A^Ty,\ \ b^Ty<0
\end{align*}
is solvable for each $s\in\{\pm1\}^n$ (is has solution $y\coloneqq y^*$). By the Farkas lemma, this is equivalent to the infeasibility of the system
\begin{align*}
(BD+A\diag(s))x=b,\ \ x\geq0.
\end{align*}
for each $s\in\{\pm1\}^n$. By substitution $x\equiv \diag(s)x$, we rewrite the system as
\begin{align*}
(A+BD\diag(s))x=b,\ \ \diag(s)x\geq0.
\end{align*}
This means that GAVE in the form $Ax+BD|x|=b$ is infeasible; see \eqref{gaveOrth}.
\end{proof}

Notice the the proof of condition \eqref{itThmGaveSolvVar3} is constructive -- it provides us with particular $D$ and $b$ such that $Ax+BD|x|=b$ is infeasible.

\begin{example}
In condition \eqref{itThmGaveSolvVar3}, the role of $D$ is unavoidable. To see it, consider, for example,
\begin{align*}
A=\begin{pmatrix}1&0\\0&1\end{pmatrix},\ \ 
B=\begin{pmatrix}1&-1\\1&0\end{pmatrix},\ \ 
D=\begin{pmatrix}1&0\\0&-1\end{pmatrix}.
\end{align*}
Then the system $Ax+B|x|=b$ does not have a unique solution for each $b\in\R^2$. However, it has at least one solution for each $b\in\R^2$. In contrast, the system $Ax+BD|x|=b$ does not have at least one solution for each $b\in\R^2$; it is infeasible for any $b<0$.
\end{example}

%%%%%
\paragraph{Special case with $A=I_n$.}
Several authors considered special GAVE in the form $x+B|x|=b$. All the results of this section apply to this form; however, we can say a bit more for this special case. 
The following conditions come from  \cite{Neu1990,Rad2016}.

\begin{theorem}\label{thmGaveSpecUnSolv}
The following conditions are equivalent:
\begin{enumerate}[(i)]
\item 
The system $x+B|x|=b$ has a unique solution for each $b\in\R^n$;
\item\label{itThmGaveSpecUnSolv1}%dle Neu1990, p. 220, pry od Rohna 
the system $x+B|x|=b$ has at most one  solution for each $b\in\R^n$;
\item\label{itThmGaveSpecUnSolv2}%dle Neu1990, p. 220, pry od Rohna 
the system $x+BD|x|=b$ has at least one  solution for each $b\in\R^n$ and each $D:|D|=I_n$;
\item\label{itThmGaveSpecUnSolv3}%dle \ref{itCorIntMcDouble2}
the system $|x|\leq|Bx|$ has only the trivial solution $x=0$;
\item\label{itThmGaveSpecUnSolv4}%dle Rad2016 pry Rohn, dle Neu1990, p. 218, pry od Rohna 
$\det(I_n+B\diag(s))>0$ for each $s\in\{\pm1\}^n$;
\item\label{itThmGaveSpecUnSolv5}%dle Rad2016 pry Rohn, dle Neu1990, p. 218, pry od Rohna 
for each $s\in\{\pm1\}^n$ and each real eigenvalue $\lambda$ of $\diag(s)B$ we have $|\lambda|<1$;
\item\label{itThmGaveSpecUnSolv6}%dle Neu1990, p. 218, pry od Rohna 
for each $s\in\{\pm1\}^n$, the inverse $(I_n+B\diag(s))^{-1}$ exists and its diagonal entries are greater than~$\frac{1}{2}$;
\item\label{itThmGaveSpecUnSolv7}%Rad2016
the map $x\mapsto x+B|x|$ is a bijection.
\end{enumerate}
\end{theorem}

\begin{proof}
Conditions \eqref{itThmGaveSpecUnSolv1}--\eqref{itThmGaveSpecUnSolv6} were presented in Neumaier~\cite[Sect.~6.1]{Neu1990}.% and attributed to Rohn.

Condition \eqref{itThmGaveSpecUnSolv7} comes from Radons~\cite{Rad2016}.
\end{proof}

From Theorem~\ref{thmGaveSpecUnSolv}\eqref{itThmGaveSpecUnSolv3} we obtain that the transformation $B\mapsto DB$ with $D\in[-I_n,I_n]$ preserves unique solvability; cf.\ Neumaier~\cite[Sect.~6.1]{Neu1990}.

\begin{corollary}\label{corGaveAspec}
The following conditions are equivalent:
\begin{enumerate}[(i)]
\item 
The system $x+B|x|=b$ has a unique solution for each $b\in\R^n$;
\item
for some $D:|D|=I_n$, the system $x+DB|x|=b$ has a unique solution for each $b\in\R^n$;
\item
for every $D:|D|=I_n$, the system $x+DB|x|=b$ has a unique solution for each $b\in\R^n$;
\item%je i v Neu1990, p. 218
for every $D\in[-I_n,I_n]$, the system $x+DB|x|=b$ has a unique solution for each $b\in\R^n$;
\item%je i v Neu1990, p. 218
the system $x+B^T|x|=b$ has a unique solution for each $b\in\R^n$.
\end{enumerate}
\end{corollary}

%%%%%%%%%%%%%%%%%%%%%%%%%%%%%%%%%%%%%%%%%%%%%%%%%%%%%%%%%%%%%%%
\section{Generalized systems}

In this section, our focus lies on exploring extensions of (generalized) absolute value equations. We commence by examining new extensions of absolute value equations, followed by an exploration of matrix absolute value equations. Subsequently, we delve into the broader realm of tensor absolute value equations, which significantly expands the scope and applicability of these equations.

%%%%%%%%%%%%%%%%%%%%%%%%%%%%%%%%%%%%%%%%%%%%%%%%%%%%%%%%%%%%%%%
\subsection{Extensions of GAVE}
%\paragraph{Extensions.}
Some further extensions of \eqref{gave} were discussed in literature \cite{KumDee2023u,Wu2021,YanWu2023,ZhoWu2021}. They are often of the type
\begin{align}\label{ggave}
Ax+B|Cx-d| = b,
\end{align}
where $A\in\R^{n\times n}$, $B\in\R^{n\times k}$, $C\in\R^{k\times n}$, $b\in\R^n$ and $d\in\R^k$. By the substitution $y\equiv Cx-d$, we can transform this system into the standard form of GAVE
\begin{align*}
Ax+B|y| = b,\ Cx-y=d.
\end{align*}
In matrix form, it reads as
\begin{align}\label{ggaveAsGave}
\begin{pmatrix}A&0\\C&-I_k\end{pmatrix}
\begin{pmatrix}x\\y\end{pmatrix}
+\begin{pmatrix}0&B\\0&0\end{pmatrix}
\begin{pmatrix}|x|\\|y|\end{pmatrix}
=\begin{pmatrix}b\\d\end{pmatrix}.
\end{align}
Thus a characterization of unique solvability of \eqref{ggave} is easily derived from those conditions on GAVE. We present some of them, those that can be a bit simplified.

\begin{corollary}\label{corIntMcTriple}
The following conditions are equivalent:
\begin{enumerate}[(i)]
\item 
The system \eqref{ggave} has a unique solution for each $b\in\R^n$ and $d\in\R^k$;
\item\label{itCorIntMcTriple0}
%for $u\in\R^n$ and $v\in\R^k$, 
the system $A^Tu+C^Tv=0$, $|v|\leq |B^Tu|$ has only the trivial solution $u=0$, $v=0$;
\item\label{itCorIntMcTriple1}
interval matrix 
\begin{align}\label{intMcTriple}
\begin{pmatrix}A&0&B\\C&-I_k&0\\0&[-I_k,I_k]&I_n\end{pmatrix}
\end{align}
is regular;
\item\label{itCorIntMcTriple2}
matrix $A + BDC$ is nonsingular for every $D\in[-I_k,I_k]$;
\item\label{itCorIntMcTriple3}
the system $Ax+By=0$, $|y|\leq|Cx|$ has only the trivial solution $x=0$, $y=0$;
\item\label{itCorIntMcTriple4}
$\det(A+B\diag(s)C)$ is constantly positive or negative for each $s\in\{\pm1\}^k$;
\item\label{itCorIntMcTriple5}
each matrix in the form 
$A+tB\diag(s)C$
is nonsingular, where $s\in\{\pm1\}^k$ and $t\in[0,1]$;
\item\label{itCorIntMcTriple6}
%each matrix in the form $A+\diag(z)$ is nonsingular, where $z_i\in[-1,1]$ for some $i\in\{1,\ldots,n\}$ and $|z_j|=1$ for $j\not=i$.
each matrix in the form 
$$A+B\diag(z^{(i)})C,\quad i=1,\ldots,k,$$
is nonsingular, where $z^{(i)}_i\in[-1,1]$ and $|z^{(i)}_j|=1$ for $j\not=i$.
\end{enumerate}
\end{corollary}

\begin{proof}
\quo{Condition \eqref{itCorIntMcTriple0}.}
By Theorem~\ref{thmGaveCharT}, the system $|A^Tu+C^Tv|\leq0$, $|v|\leq |B^Tu|$ should have only the trivial solution.

\quo{Condition \eqref{itCorIntMcTriple1}.}
By Theorem~\ref{thmGaveRegABD}, the interval matrix
\begin{align*}
\begin{pmatrix}A&0&0&B\\C&-I_k&0&0\\{}[-I_n,I_n]&0&I_n&0\\0&[-I_k,I_k]&0&I_k\end{pmatrix}
\end{align*}
is regular. By permutation, this matrix takes the form
\begin{align*}
\begin{pmatrix}A&0&B&0\\C&-I_k&0&0\\0&[-I_k,I_k]&I_k&0\\{}[-I_n,I_n]&0&0&I_n\end{pmatrix},
\end{align*}
which is block lower triangular. So the first $3\times 3$ diagonal matrix block must be regular; the second diagonal block is simply $I_n$.

\quo{Condition \eqref{itCorIntMcTriple2}.}
Let $D\in[-I_k,I_k]$ be arbitrary. By elementary operations, we obtain
\begin{align*}
\begin{pmatrix}A&0&B\\C&-I_k&0\\0&D&I_k\end{pmatrix}
\sim\begin{pmatrix}A&-BD&0\\C&-I_k&0\\0&D&I_k\end{pmatrix}
\end{align*}
By column elementary operations on the first block, we obtain
\begin{align*}
\begin{pmatrix}A&-BD\\C&-I_k\end{pmatrix}
\sim\begin{pmatrix}A-BDC&-BD\\0&-I_k\end{pmatrix}.
\end{align*}

\quo{Condition \eqref{itCorIntMcTriple3}.}
%By Corollary~\ref{corIntMcDouble}\ref{itCorIntMcDouble2} 
By Rohn \cite[cond.~(ii)]{Roh2009} applied on \eqref{intMcTriple} we get the system
$$
|Ax+By|\leq0,\ |Cx-z|\leq0,\ |y|\leq|z|,
$$
from which we have the result by the substitution $z\equiv Cx$.

\quo{Condition \eqref{itCorIntMcTriple4}.}
This follows from Rohn \cite[cond.~(xxxii)]{Roh2009} .

\quo{Condition \eqref{itCorIntMcTriple5}.}
This follows from Rohn \cite[cond.~(xxxviii)]{Roh2009}.

\quo{Condition \eqref{itCorIntMcTriple6}.}
This follows from Rohn \cite[cond.~(xli)]{Roh2009}.
\end{proof}

It is an open question whether the unique solvability of \eqref{ggave} can also be deduced from the special case of $d=0$. That is, whether unique solvability of $Ax+B|Cx| = b$ for each $b\in\R^n$ implies unique solvability of \eqref{ggave}. Anyway, one step towards this property is formulated in the following statement.

\begin{theorem}\label{thmGgaveAtLeast}
The following conditions are equivalent:
\begin{enumerate}[(i)]
\item\label{itThmGgaveAtLeast1}
The system \eqref{ggave} has a unique solution for each $b\in\R^n$ and $d\in\R^k$;
\item\label{itThmGgaveAtLeast2}
the system $Ax+BD|Cx|=b$ has at least one solution for each $b\in\R^n$ and each $D:|D|=I_k$.
\end{enumerate}
\end{theorem}

\begin{proof}
\quo{\eqref{itThmGgaveAtLeast1}$\Rightarrow$\eqref{itThmGgaveAtLeast2}.} 
This implication is clear in view of Theorem~\ref{thmGaveSolvVar}\eqref{itThmGaveSolvVar3}.

\quo{\eqref{itThmGgaveAtLeast2}$\Rightarrow$\eqref{itThmGgaveAtLeast1}.} 
We prove it by contraposition. If condition \eqref{itThmGgaveAtLeast1} is not satisfied, then by Corollary~\ref{corIntMcTriple}\eqref{itCorIntMcTriple0} there exist nontrivial $u^*$ and $v^*$ such that $A^Tu^*+C^Tv^*=0$, $|v^*|\leq |B^Tu^*|$. Necessarily we have $u^*\not=0$, otherwise $v^*=0$. Now, in a similar fashion as in the proof of Theorem~\ref{thmGaveSolvVar}\eqref{itThmGaveSolvVar3}, we define $b\coloneqq -u^*$, $d\coloneqq0$ and $D\coloneqq\diag(\sgn(B^Tu^*))$. Then the system
\begin{align*}
A^Tu+C^Tv=0,\ \ 
DB^Tu \geq \diag(s)v,\ \ b^Tu<0
\end{align*}
is solvable for each $s\in\{\pm1\}^k$ (is has solution $(u,v)\coloneqq (u^*,v^*)$). By the Farkas lemma, this is equivalent to the infeasibility of the system
\begin{align*}
Ax+BDy=b,\ \ 
Cx-\diag(s))y=0,\ \ y\geq0.
\end{align*}
for each $s\in\{\pm1\}^k$. By substitution $y\equiv \diag(s)Cx$, we rewrite the system as
\begin{align*}
Ax+BD\diag(s)Cx=b,\ \ \diag(s)Cx\geq0,
\end{align*}
which is equivalent to infeasibility of $Ax+BD|Cx|=b$.
\end{proof}

%%%%%
\paragraph{Special case with $B=I_n$.}
Particularly for the system $Ax-|Cx|=b$ with $A,~C\in\R^{n \times n}$, Wu~\cite{Wu2021} proposed the following characterization. For the sake of uniformity of presentation, we formulate it for the form $Ax+|Cx|=b$.

\begin{comment}
\begin{theorem}
The following conditions are equivalent:
\begin{enumerate}[(i)]
\item 
The system $Ax-|Bx|=b$ has a unique solution for each $b\in\R^n$;
\item
$\{A-B,A+B\}$ has the row $\mathcal{W}$-property;
\item
$A+B$ is nonsingular and $\{I_n,(A-B)(A+B)^{-1}\}$ has the row $\mathcal{W}$-property.
\end{enumerate}
\end{theorem}
\end{comment}

\begin{theorem}\label{thmGaveModif}
The following conditions are equivalent:
\begin{enumerate}[(i)]
\item 
The system $Ax+|Cx|=b$ has a unique solution for each $b\in\R^n$;
\item
$\{A-C,A+C\}$ has the row $\mathcal{W}$-property;
\item
$A-C$ is nonsingular and $\{I_n,(A+C)(A-C)^{-1}\}$ has the row $\mathcal{W}$-property;
\item
$A-C$ is nonsingular and $(A+C)(A-C)^{-1}$ is a P-matrix.
\end{enumerate}
\end{theorem}

Kumar and Deepmala \cite{KumDee2023u} derived the same conditions for unique solvability of the system $Ax+|Cx-d|=b$ with $A,~C\in\R^{n \times n}$. This justifies the following statement, which answers the open problem stated before Theorem~\ref{thmGgaveAtLeast} for the case when $B=I_n$.

\begin{theorem}
The following conditions are equivalent:
\begin{enumerate}[(i)]
\item 
The system $Ax+|Cx|=b$ has a unique solution for each $b\in\R^n$;
\item
the system $Ax+|Cx-d|=b$ has a unique solution for each $b,d\in\R^n$.
\end{enumerate}
\end{theorem}

%%%%%
\paragraph*{Further variants.}
Another form of absolute value systems investigated are (generalized) absolute value equations associated with second order cones. Visually, they have the same form as GAVE, that is, $Ax+B|x|=b$. However, the difference is in the meaning of the absolute value. Herein, it is defined as $|x|=\sqrt{x\circ x}$, where $\circ$ denotes the Jordan product associated with the Lorentz cone. For such systems, there are known sufficient conditions for unique solvability \cite{HuHua2011,HuaLi2022,MiaHsu2021}, however, a complete characterization is an open problem.

%%%%%%%%%%%%%%%%%%%%%%%%%%%%%%%%%%%%%%%%%%%%%%%%%%%%%%%%%%%%%%%
\subsection{Matrix absolute value equations}\label{sMave}

Matrix absolute value equations have the form $AX+|X|=F$, and the generalized version reads as $AX+B|X|=F$. It will not make the description and derivation complicated if we consider the more general form \cite{Has2021,MolBei2022,TanMia2022,WanLi2021}
\begin{align}\label{mave}
AXC+B|X|E = F,
\end{align}
where $A,B\in\R^{m\times n}$, $C,E\in\R^{p\times q}$ and $F\in\R^{m\times q}$ are input matrices and $X\in\R^{n\times p}$ is unknown. Characterization of unique solvability then follows by the transformation of the matrix equation system into the form of GAVE 
\begin{align}\label{maveAsGave}
(C^T\otimes A)\vek(X)+(E^T\otimes B)\vek(|X|)=\vek(F).
\end{align}
Herein, $\otimes$ denotes the Kronecker product and $\vek(\cdot)$ is the vectorization operator that converts the matrix into a vector,
$$
\vek(F)=(F_{11},\dots,F_{m1},F_{12},\dots,F_{m2},\dots,F_{1q},\dots,F_{mq})^T;
$$
see, e.g., Horn and Johnson \cite{HorJoh1991}. Such a transformation was utilized by Wang \& Li~\cite{WanLi2021} or Tang \& Miao~\cite{TanMia2022} to derive certain necessary and sufficient conditions for unique solvability. In principle, we can apply all the conditions for GAVE from Section~\ref{sGave} to characterize the unique solvability of~\eqref{maveAsGave}. Nevertheless, the size of the system \eqref{maveAsGave} grows up to $mq\times mq$, so we will be particularly interested in conditions that work with the initial size and avoid the Kronecker product.

\begin{theorem}\label{thmMaveCharac}
The following conditions are equivalent:
\begin{enumerate}[(i)]
\item 
The system \eqref{mave} has a unique solution for each $F\in\R^{m\times q}$;
\item\label{itThmMaveCharac2}
the system $|A^TXC^T|\leq|B^TXE^T|$ has only the trivial solution $X=0$;
\item\label{itThmMaveCharac3}
the system $AXC+BYE=0$, $|Y|\leq|X|$ has only the trivial solution $X=Y=0$.
\end{enumerate}
\end{theorem}

\begin{proof}
\quo{Condition \eqref{itThmMaveCharac2}.}
By Theorem~\ref{thmGaveCharT} and using the form \eqref{maveAsGave}, we obtain the system
\begin{align*}
|(C^T\otimes A)^T\vek(X)| \leq |(E^T\otimes B)^T\vek(X)|,
\end{align*}
 or
\begin{align*}
|(C\otimes A^T)\vek(X)| \leq |(E\otimes B^T)\vek(X)|.
\end{align*}
Its matrix version is $|A^TXC^T|\leq|B^TXE^T|$ with variable $X\in\R^{m\times q}$.

\quo{Condition \eqref{itThmMaveCharac3}.}
It follows from Corollary~\ref{corIntMcDouble}\eqref{itCorIntMcDouble2}. The system $(C^T\otimes A)\vek(X)+(E^T\otimes B)\vek(Y)=0$, $|\vek(Y)|\leq|\vek(X)|$ is simply reformulated in the matrix form. 
\end{proof}

The other conditions seem to be more difficult to adapt for the matrix equations. However, for the special type in the form $AX+B|X|=F$ with $A,B,F\in\R^{n\times n}$, we can derive stronger results. The theorem below states that all conditions from Section~\ref{sGave} apply with no modification.

\begin{theorem}\label{thmMaveCharacGave}
The following conditions are equivalent:
\begin{enumerate}[(i)]
\item 
The system $AX+B|X|=F$ has a unique solution for each $F\in\R^{n\times n}$;
\item\label{itThmMaveCharacGave1}
matrix $A + BD$ is nonsingular for every $D\in[-I_n,I_n]$;
\item\label{itThmMaveCharacGave2}
the GAVE system $Ax+B|x|=b$ has a unique solution for each $b\in\R^n$;
\item\label{itThmMaveCharacGave3}
for any upper triangular $C,E$ with ones on the diagonals, the system $AXC+B|X|E = F$ has a unique solution for each $F\in\R^{n\times n}$.
\end{enumerate}
\end{theorem}

\begin{proof}
\quo{Condition \eqref{itThmMaveCharacGave1}.}
By Theorem~\ref{thmGaveUnSol}, for every $D\in[-I_{n^2},I_{n^2}]$ the matrix $(I_n\otimes A) + (I_n\otimes B)D$ is nonsingular. This master matrix is block diagonal, whence each block matrix $A+BD_i$ should be nonsingular for every $D_i\in[-I_n,I_n]$, $i=1,\dots,n$. Thus it is sufficient to consider just one block.

\quo{Condition \eqref{itThmMaveCharacGave2}.}
The previous condition~\eqref{itThmMaveCharacGave1} is the same as that in Theorem~\ref{thmGaveUnSol}.

\quo{Condition \eqref{itThmMaveCharacGave3}.}
Analogous to condition \eqref{itThmMaveCharacGave1}.
\end{proof}

Slightly different system was considered by Li~\cite{Li2022}, who addressed the matrix system of type $AXC+|BXE| = F$ with $A,B\in\R^{m\times n}$, $C,E\in\R^{p\times q}$ and $F\in\R^{m\times q}$. The proposed unique solvability conditions are analogous to those from Theorem~\ref{thmGaveModif} applied on the master system
\begin{align*}
(C^T\otimes A)\vek(X)+|(E^T\otimes B)\vek(X)|=\vek(F).
\end{align*}

%%%%%
\paragraph*{Extensions.}
The characterization of unique solvability of the absolute value matrix system \eqref{mave} can directly be extended to cover a more general system in the form of
\begin{align}\label{gmave}
\sum_{i=1}^{i^*} A_iXC_i+\sum_{j=1}^{j^*} B_j|X|E_j = F,
\end{align}
where $A_i,B_j\in\R^{m\times n}$, $C_i,E_j\in\R^{p\times q}$ and $F\in\R^{m\times q}$. We present the conditions without the proof, which is straightforward.

\begin{theorem}
The following conditions are equivalent:
\begin{enumerate}[(i)]
\item 
The system \eqref{gmave} has a unique solution for each $F\in\R^{n\times n}$;
\item
the system $|\sum_{i=1}^{i^*}A^T_iXC^T_i|\leq|\sum_{j=1}^{j^*}B^T_jXE^T_j|$ has only the trivial solution $X=0$;
\item
the system $\sum_{i=1}^{i^*} A_iXC_i+\sum_{j=1}^{j^*} B_jYE_j=0$, $|Y|\leq|X|$ has only the trivial solution $X=Y=0$.
\end{enumerate}
\end{theorem}

%%%%%
\paragraph*{Tensor absolute value equations.}

In recent studies, there has been a growing interest in investigating tensor absolute value equations (TAVEs) as a further generalization. Various researchers have proposed sufficient conditions and special tensor classes to address the issue of unique solvability in these equations \cite{BeiKal2022, CuiLia2022, DuZha2018, JiaLi2021}. 
For instance, in a related paper \cite{cui2022existence}, the authors investigate the existence and uniqueness of solutions for both tensor complementarity problems (TCP) and TAVEs with a special structure. They  provide sufficient conditions for the existence and uniqueness of solutions for TAVEs with a special structure.
Regarding the unique solvability of TAVEs, although various sufficient conditions and special tensor classes have been proposed, a complete characterization of the necessary and sufficient conditions is still lacking. Therefore, further research is needed to achieve a comprehensive understanding of the unique solvability of tensor absolute value equations.

%%%%%%%%%%%%%%%%%%%%%%%%%%%%%%%%%%%%%%%%%%%%%%%%%%%%%%%%%%%%%%%
\section{Conclusion and future directions}

This paper surveys necessary and sufficient conditions for the unique solvability of absolute value equations and their several generalizations.  We not only present the known conditions but also introduce some novel ones that have not been published before. We also tried to find the tight frontiers -- by counter examples we showed that some conditions valid for special cases need not be satisfied for their generalizations. 

The results presented in this paper are not solely of theoretical interest. Some of the conditions have already been used to derive cheap sufficient solvability conditions. So possibly the others can also serve as a basis to develop other tractable sufficient or necessary conditions.

 However, there is still room for further exploration and the formulation of alternative equivalent characterizations. Several open problems remain, including the one mentioned prior to Theorem~\ref{thmGgaveAtLeast}, as well as the complete characterization of the unique solvability of tensor absolute value equations. The development of robust solution methodologies under uncertain data and the exploration of scalable implementations are among the avenues of future investigation. By addressing these open problems, we can enrich the understanding and applicability of these mathematical models. One key aspect to consider is the characterization of the set of feasible solutions for uncertain AVE and GAVE and analyzing the computational complexity of solving robust AVE.

%%%%%%%
\paragraph{Acknowledgments.} 
 The research work of Shubham Kumar was supported by the Ministry of Education, Government of India, through Graduate Aptitude Test in Engineering (GATE) fellowship registration No. MA19S43033021. The work of Milan Hlad\'{i}k was supported by the Czech Science Foundation Grant P403-22-11117S.

%%%%%%%%%%%%%%%%%%%%%%%%%%%%%%%%%%%%%%%%%%%%%%%%%%%%%%%%%%%%%%%
% REFERENCES
%%%%%%%%%%%%%%%%%%%%%%%%%%%%%%%%%%%%%%%%%%%%%%%%%%%%%%%%%%%%%%%

\bibliographystyle{abbrv}
\bibliography{survey_ave_solv}

\end{document}